%% file: RSBarxiv.tex
\documentclass[10pt,twoside]{article}
\usepackage{amsmath}
\usepackage{amsfonts}
\usepackage{amssymb}
\usepackage{vector}
\def\eps{\varepsilon}
\font\tencmmib=cmmib10 \skewchar\tencmmib '60
\newfam\cmmibfam
\textfont\cmmibfam=\tencmmib

\def\bbox{\quad\hbox{\vrule \vbox{\hrule \vskip2pt \hbox{\hskip2pt
\vbox{\hsize=1pt}\hskip2pt} \vskip2pt\hrule}\vrule}}
\def\lessim{\ \lower4pt\hbox{$
\buildrel{\displaystyle <}\over\sim$}\ }
\def\gessim{\ \lower4pt\hbox{$\buildrel{\displaystyle >}
\over\sim$}\ }
\def\P{{\cal P}}

\def\M{{\cal M}}

\def\eps{{\varepsilon}}
\def\ch{{\mbox{ch}}}

\def\la{{\Bigl\langle}}
\def\ra{{\Bigr\rangle}}

\def\qed{\hfill\break\rightline{$\bbox$}}
\newcommand{\e}{\mathbb{E}}

\newcommand{\Reals}{\mathbb{R}}

\newcommand{\vsi}{{\vec{\sigma}}}

\newtheorem{lemma}{Lemma}
\newtheorem{theorem}{Theorem}

\makeatletter
\@addtoreset{equation}{section}

\makeatother

\input{invlat.tex}

\begin{document}

\title{On differentiability of the Parisi formula.}

\author{ 
Dmitry Panchenko\thanks{Department of Mathematics, Texas A\&M University,
email: panchenk@math.tamu.edu. This work is partially supported by NSF grant.
}\\
{\it Department of Mathematics}\\
{\it Texas A\&M University}\\
}
\date{}

\maketitle

\begin{abstract}
It was proved by Michel Talagrand in \cite{PM} that the Parisi formula for the free energy in the Sherrington-Kirkpatrick model is differentiable with respect to inverse temperature parameter. We present a simpler proof of this result by using approximate solutions in the Parisi formula and give one example of application of the differentiability to prove non self-averaging of the overlap outside of the replica symmetric region.
\end{abstract}
\vspace{0.5cm}



\setcounter{section}{1}
\subsection{Introduction and main results.}

Let us consider a $p$-spin Sherrington-Kirkpatrick Hamiltonian
$$
H_{N,p}(\vsi)=\frac{1}{N^{(p-1)/2}}\sum_{1\leq i_1,\ldots,i_p \leq N}
g_{i_1,\ldots,i_p}\sigma_{i_1}\ldots \sigma_{i_p}
$$ 
indexed by spin configurations $\vsi\in\Sigma_N = \{-1,+1\}^N$ 
where $(g_{i_1,\ldots,i_p})$ are i.i.d. standard Gaussian random variables.
A mixed $p$-spin Hamiltonian is defined as the sum
\begin{equation}
H_N(\vsi) = \sum_{p\geq 1} \beta_p\, H_{N, p}(\vsi)
\label{HN}
\end{equation}
over a finite set of indices $p\geq 1.$ 
The covariance of $H_N$ can be easily computed
\begin{equation}
\e H_N(\vsi^1) H_N(\vsi^2) = N \xi(R_{1,2}),
\label{cov}
\end{equation}
where 
$$
R_{1,2}=\frac{1}{N}\sum_{i\leq N}\sigma_i^1 \sigma_i^2
\,\,\,\mbox{ and }\,\,\,
\xi(x)=\sum_{p\geq 1} \beta_p^2 x^p.
$$ 
A quantity $R_{1,2}$ is called the overlap of configurations $\vsi^1, \vsi^2.$ 
To avoid the trivial case when all the spins decouple we assume that $\beta_p\not = 0$ for at least one $p\geq 2$ so that
$\xi''(x)>0$ for $x > 0.$
Given an external field parameter $h\in\Reals,$ the free energy is defined by 
\begin{equation}
F_N(\vec{\beta}) = \frac{1}{N}\,\e\log \sum_{\vsi} 
\exp\bigl(H_N(\vsi)+h\sum_{i\leq N} \sigma_i\bigr).
\label{FE}
\end{equation}
The problem of computing the thermodynamic limit of the free energy $\lim_{N\to\infty} F_N$
is one of the central questions in the analysis of the SK model and the value of this limit was 
predicted by Giorgio Parisi in \cite{Parisi} as a part of his celebrated theory
that goes far beyond the computation of the free energy. 
The prediction of Parisi was confirmed with mathematical rigor by Michel Talagrand in \cite{T-P} 
following a breakthrough of Francesco Guerra in \cite{Guerra} where
a replica symmetry breaking interpolation was introduced.
Validity of the Parisi formula provides a lot of information about the model and,
in particular, about the distribution of the overlap under the Gibbs measure
corresponding to the Hamiltonian $H_N(\vsi)$. In the next section we will show
one important application of the Parisi formula which is based on its differentiability 
with respect to inverse temperature parameters. Namely,
we will prove a stronger version of the result of Pastur and Shcherbina in \cite{PS} about 
the non self-averaging of the overlap at low temperature. 

In the remainder of this section we present a simplified version of the argument 
of Talagrand in \cite{PM} and prove the differentiability of the Parisi formula. 
Let us start by recalling the definition of the Parisi formula.
Let $\M$ be the set of cumulative distribution functions on $[0,1].$
We will identify a c.d.f. $m$ with a distribution it defines and simply 
call $m$ itself a distribution on $[0,1].$ 
A distribution with at most $k$ atoms is defined by
\begin{equation}
m(q)= \sum_{0\leq l\leq k} m_l I( q_l \leq q < q_{l+1})
\label{mk}
\end{equation}
for some sequences
$$
0=m_0\leq m_1\leq \ldots\leq m_{k-1}\leq m_k=1,
$$
$$
0=q_0\leq q_1\leq \ldots\leq q_k\leq q_{k+1} = 1.
$$
Consider independent Gaussian r.v. $(z_l)_{0\leq l \leq k}$
such that $\e z_{l}^2 = \xi'(q_{l+1}) - \xi'(q_l).$
Let
$$
X_{k} = \log \ch\Bigl(\sum_{0\leq l\leq k}z_l + h\Bigr)
$$
and recursively for $1\leq l\leq k$ define
\begin{equation}
X_{l-1} = \frac{1}{m_l}\log \e_l \exp m_l X_{l}
\label{recX}
\end{equation}
where $\e_l$ denotes the expectation in $(z_{p})$ for $l\leq p\leq k.$
Define
\begin{equation}
\P(m,\vec{\beta}) = \e X_0
-\frac{1}{2}\sum_{1\leq l\leq k}m_l(\theta(q_{l+1}) - \theta(q_l)).
\label{Parisi}
\end{equation}
where $\theta(x)=x\xi'(x)-\xi(x).$ 
On the set of discrete $m\in\M$ as in (\ref{mk}) 
the functional $\P(m,\vec{\beta})$ is Lipschitz in $m$ 
with respect to $L_1$ norm (see \cite{Guerra}, \cite{PM}).
Therefore, it can be extended by continuity to a Lipschitz
functional on the entire space $\M.$
The Parisi formula is then defined by
\begin{equation}
\P(\vec{\beta})=\inf_{m\in \M}\P(m,\vec{\beta}).
\label{Par}
\end{equation}
This infimum is obviously achieved by continuity and compactness.
Any $m\in\M$ that achieves the infimum is called a Parisi measure.
It is conjectured (\cite{Parconv}) that $\P(m,\vec{\beta})$ is convex in $m$  in which
case the Parisi measure would be unique. 

By H\"older's inequality, $F_N(\vec{\beta})$ is convex in $\vec{\beta}$ 
and, thus, its limit $\P(\vec{\beta})$ is also convex. Convexity implies that
$\P(\vec{\beta})$ is differentiable in each parameter $\beta_p$ almost everywhere and 
it was proved in \cite{PM} that $\P(\vec{\beta})$ is in fact differentiable
for all values of $\beta_p\,$. 
The proof was based on a careful analysis of 
the functional $\P(m,\vec{\beta})$ in the neighborhood of a Parisi measure
and parts of the proof were rather technical due to the fact that a Parisi measure
is not necessarily discrete. 
We will prove a slightly weaker analogue of Theorem 1.2 in \cite{PM} but we will bypass these difficulties 
by working with approximations of a Parisi measure by discrete measures of the type (\ref{mk}). 
The main difference is that we express the derivative in (\ref{derP}) below in terms of some Parisi 
measure instead of any Parisi measure as in \cite{PM}.

\begin{theorem}\label{Th1}
The derivative of the Parisi formula $\P(\vec{\beta})$ with respect to any $\beta_p$ exists and
\begin{equation}
\frac{\partial \P(\vec{\beta})}{\partial \beta_p} = 
\beta_p\Bigl(1-\int q^p dm_{\vec{\beta}}(q)\Bigr)
\,\,\,\mbox{ for all }\,\,\, p\geq 1
\label{derP}
\end{equation}
for some Parisi measure $m_{\vec{\beta}}$.
\end{theorem}

To prove Theorem \ref{Th1} we will first obtain a similar statement for discrete approximations 
of a Parisi measure; this result corresponds to Proposition 3.2 in \cite{PM}.

\begin{lemma}\label{Lem1}
Given $k\geq 1,$ suppose that $m\in \M$ achieves the minimum of $\P(m,\vec{\beta})$ over 
all distributions with at most $k$ atoms as in (\ref{mk}). Then
$$
\frac{\partial\P}{\partial \beta_p}(m,\vec{\beta})
=
\beta_p\Bigl(1-\int q^{p} dm(q)\Bigr).
$$
\end{lemma}
\textbf{Proof.}
Suppose that $m$ has $k'$ atoms in $(0,1)$ for some $k'\leq k.$
For simplicity of notations, let us assume that $k'=k.$
Let us start by noting that $\e X_0$ depends on $\vec{\beta}$ only through 
$\xi'(1)$ and $\xi'(q_l)$ for $1\leq l\leq k.$  Let us make
the dependence on $\xi'(1)$ explicit. Since 
$$
X_{k-1} = \log \ch\bigl(\sum_{0\leq l\leq k-1} z_l + h\bigr)
+\frac{1}{2}(\xi'(1)-\xi'(q_k))
$$
we can continue recursive construction (\ref{recX}) to show that
$$
\e X_0 = \frac{1}{2}\,\xi'(1) + \frac{1}{2}\,f(\xi'(q_1),\ldots,\xi'(q_k))
$$
for some smooth function $f(x_1,\ldots,x_k):\Reals^k\to\Reals.$
Then, rearranging the terms in (\ref{Parisi})
\begin{equation}
\P(m,\vec{\beta})
=
\frac{1}{2}\,\xi(1) + \frac{1}{2}\,f(\xi'(q_1),\ldots,\xi'(q_k))
+\frac{1}{2}\sum_{1\leq l\leq k}(m_l-m_{l-1})\theta(q_l). 
\label{Pmb}
\end{equation}
Since $m$ achieves the minimum, for $1\leq l\leq k$
$$
2\frac{\partial\P}{\partial q_l}
=
\frac{\partial f}{\partial x_l}\, \xi''(q_l) +(m_l - m_{l-1})q_l\xi''(q_l)=0
$$
and since $\xi''(q)>0$ for $q>0$ this implies that
\begin{equation}
\frac{\partial f}{\partial x_l} 
= -(m_l - m_{l-1})q_l.
\label{Phix}
\end{equation}
Since 
$$
\xi(q)=\sum_{p\geq 1} \beta_p^2\, q^p,\,\,\,
\xi'(q)=\sum_{p\geq 1} p\,\beta_p^2\, q^{p-1} 
\,\,\mbox{ and }\,\,
\theta(q)=\sum_{p\geq 1} (p-1)\beta_p^2\, q^p,
$$ 
using (\ref{Pmb}) and (\ref{Phix}) we compute
\begin{eqnarray*}
\frac{\partial\P}{\partial \beta_p}
&=&
\beta_p
+ \sum_{1\leq l\leq k}
\frac{\partial f}{\partial x_l}\, p\,\beta_p\, q_l^{p-1} 
+\sum_{1\leq l \leq k}(m_l - m_{l-1})(p-1)\beta_p\, q_l^p
\\
&=&
\beta_p - \beta_p \sum_{1\leq l \leq k}(m_l - m_{l-1})q_l^p
=
\beta_p\Bigl(1-\int q^{p} dm(q)\Bigr)
\end{eqnarray*}
and this finishes the proof.
\qed

\textbf{Proof of Theorem \ref{Th1}.} 
First of all, let us fix all but one parameter in $\vec{\beta}$ 
and think of all the functions that depend on $\vec{\beta}$ as functions
of one variable $\beta=\beta_p.$
Let $m^k$ be a distribution from Lemma \ref{Lem1}. By definition of Parisi formula
and Lipschitz property of $\P(m,\beta)$ we have
$\P(m^k,\beta)\downarrow \P(\beta)$ as $k\to\infty$ or, in other words,
\begin{equation}
0\leq \P(m^k,\beta) - \P(\beta)\leq \eps_k
\label{PNP}
\end{equation}
for some sequence $\eps_k\downarrow 0.$
To prove that a convex function $\P(\beta)$ is differentiable we need
to show that its subdifferential $\partial\P(\beta)$ contains a unique point. 
Let $a\in\partial \P(\beta).$ Then by convexity of $\P$, (\ref{PNP}) and
the fact that $\P(\beta')\leq \P(m^k,\beta')$ for all $\beta',$
$$
a\leq \frac{\P(\beta+y) - \P(\beta)}{y}
\leq
\frac{\P(m^k, \beta+y) - \P(m^k,\beta) + \eps_k}{y}
$$
and
$$
a\geq \frac{\P(\beta) - \P(\beta-y)}{y}
\geq
\frac{\P(m^k, \beta) - \P(m^k,\beta-y) - \eps_k}{y}
$$
for $y>0.$
It is a simple exercise to check that for any discrete $m\in \M$
the second derivative 
$\partial^2 \P(m,\beta)/\partial\beta^2$
stays bounded if $\beta$ stays bounded and the bound is uniform
in $m$ (see \cite{T-P} or \cite{PM}).
Therefore, using Taylor's expansion around $y=0$ on the right hand side of the above inequalities
gives
$$
\frac{\partial\P}{\partial\beta}(m^k,\beta) - Ly - \frac{\eps_k}{y}
\leq
a
\leq 
\frac{\partial \P}{\partial\beta}(m^k,\beta) + Ly + \frac{\eps_k}{y}.
$$
Taking $y=\sqrt{\eps_k}$ we obtain
$$
a=
\frac{\partial\P}{\partial\beta}(m^k,\beta) + {\cal O}(\sqrt{\eps_k})
=
\beta\Bigl(1-\int q^{p} dm^k(q)\Bigr)+ {\cal O}(\sqrt{\eps_k})
$$
by Lemma \ref{Lem1}. Finally, taking a subsequence of $(m^k)$ that converges 
in $L_1$ norm to some Parisi measure $m_{\vec{\beta}}$ proves that
$$
a 
= 
\beta\Bigl(1-\int q^{p} dm_{\vec{\beta}}(q)\Bigr).
$$  
This uniquely determines $a$ and, thus, $a=\P'(\beta).$
\qed

\subsection{Non self-averaging of the overlap.}

In this section we make an assumption that all indices in (\ref{HN}) 
are even numbers with one possible exception of $p = 1,$ i.e. besides
a trivial linear term we consider only even spin interaction terms. 
The reason for this is because the validity of the Parisi formula 
was proved in \cite{T-P} under certain conditions on the function $\xi$
which essentially correspond to the choice of only even spin interaction
terms. Under this assumption, by \cite{T-P}, 
$$
\lim_{N\to\infty} F_N(\vec{\beta}) = \P(\vec{\beta})
$$ 
and since both
$F_N(\vec{\beta})$ and $\P(\vec{\beta})$ are convex functions and, by Theorem \ref{Th1}, 
$\P(\vec{\beta})$ is differentiable in $\beta_p,$ we get
$$
\lim_{N\to\infty} 
\frac{\partial F_N}{\partial\beta_p}
=
\frac{\partial\P}{\partial\beta_p}
=
\beta_p\Bigl(1-\int q^{p} dm_{\vec{\beta}}(q)\Bigr).
$$
By Gaussian integration by parts one can easily see that,
$$
\frac{\partial F_N}{\partial\beta_p}
= \beta_p\Bigl(1-\e\bigl\la R_{1,2}^{p}\bigr\ra\Bigr)
$$
where $\la\cdot\ra$ is the Gibbs average with respect to the
Hamiltonian $H_{N}(\vsi)$ and, therefore, 
for any $p\geq 1$  such that $\beta_p>0$ we get
\begin{equation}
\lim_{N\to\infty }\e \la R_{1,2}^{p}\ra = \int q^{p} dm_{\vec{\beta}}(q).
\label{derapprox}
\end{equation}
Thus, from Theorem \ref{Th1} one obtains information about moments
of the overlap, in particular, about the existence of their thermodynamic limit.
(This result is not new, it appears in \cite{PFOP} and \cite{PM}.)
If Hamiltonian $H_N(\vsi)$ contains all even $p$-spin interaction terms then (\ref{derapprox})
holds for all even $p\geq 2$ and, thus, the distribution of $|R_{1,2}|$ is approximated
by the Parisi measure $m_{\vec{\beta}}.$ It is predicted by the Parisi theory that this is also 
true when only a finite number of even $p$-spin interaction terms are present;
however, this is an open problem. (\ref{derapprox}) provides information only about the moments of
the overlap corresponding to the terms present in the Hamiltonian.

We will now use this information to give two examples of non self-averaging of the overlap.
To put these examples in perspective, let us first recall several well-known results
about the classical $2$-spin SK model, $H_N=\beta H_{N,2},$ without
external field, $h=0.$  
Let us recall that inverse temperature parameter $\beta$ is said to belong to replica 
symmetric region if the infimum in the Parisi formula (\ref{Par}) is achieved 
on Dirac measure $\delta_0$ concentrated at zero.  In this simplest case
the Parisi formula $\P(\beta)$ is called a replica symmetric solution.
It was proved by Aizenman, Lebowitz and Ruelle in \cite{ARL} that replica symmetric solution 
holds for $\beta^2\leq 2$ and it was proved by Toninelli in \cite{TAT} that it does not hold 
for $\beta^2>2$ (the result in \cite{TAT} is more general, it also covers the case with external field).
In other words, the set of $\beta^2\leq 2$ is the replica symmetric region.
Note that the reason we have $\beta^2\leq 2$ instead of a more familiar $\beta^2\leq 1$ 
is because for simplicity  we defined the Hamiltonian $H_{N,2}$ as the sum over all indices 
$i_1$ and $i_2$ rather than $i_1<i_2.$
A well-known result of Pastur and Shcherbina in \cite{PS} states that if 
\begin{equation}
\lim_{N\to\infty}\e(\la R_{1,2}\ra - \e\la R_{1,2}\ra)^2 = 0
\label{PastS}
\end{equation}
then replica symmetric solution holds. Therefore, for $\beta^2>2$ (\ref{PastS}) can not hold
and this implies that $\limsup_{N\to\infty}\e \la R_{1,2}^2 \ra >0.$ 
Differentiability of the Parisi formula implies that the limit $\lim_{N\to\infty}\e \la R_{1,2}^2 \ra$ 
in (\ref{derapprox}) exists and, consequently, the result of Pastur
and Shcherbina can be used to deduce that this limit is strictly positive when $\beta^2>2$.
However, one can give a more direct proof of a more general result without invoking \cite{PS}.

\textbf{Example 1} ($h=0, \beta_1=0$). This case is similar to the classical SK model
without external field, only now $p$-spin interactions for even $p>2$ are also allowed.
A replica symmetric region is again defined as the set of parameters
$\vec{\beta}$ such that the infimum in (\ref{Par}) is achieved on Dirac
measure $\delta_0$ concentrated at zero, but the description of this region is slightly
more complicated (see Theorem 2.11.16 in \cite{SG}).
Using the continuity of the functional $m\to\P(m,\vec{\beta})$
with respect to the $L_1$ norm (see \cite{Guerra}, \cite{PM}),
outside of the replica symmetric region any Parisi measure $m_{\vec{\beta}}$ must satisfy
$m_{\vec{\beta}}(\{q >0\})>0$.
Therefore, by (\ref{derapprox}), for any even $p\geq 2$ such that $\beta_p>0$ we have  
\begin{equation}
\lim_{N\to\infty }\e \la R_{1,2}^{p}\ra > 0.
\end{equation}
Since by symmetry, $\la R_{1,2}\ra = 0,$ this proves non self-averaging of the overlap
outside of the replica symmetric region.

\textbf{Example 2}  ($h\not = 0,$ $\beta_{p_1},\beta_{p_2}\not = 0$ for some $p_1<p_2$).
A similar argument can be used in the
presence of external field if at least two different even $p$-spin interaction
terms are present. In this case, due to the absence of symmetry, a replica symmetric region is defined 
as the set of parameters $\vec{\beta}$ such that the infimum in (\ref{Par}) is achieved on Dirac
measure $\delta_x$ concentrated at any point $x\in [0,1]$ rather than zero.
Again, by continuity of $m\to \P(m,\vec{\beta}),$
on the complement of the replica symmetric region any Parisi measure $m_{\vec{\beta}}$ must satisfy
$$
\int |q - x|dm_{\vec{\beta}}(q)\geq \eps 
$$
for all $x\in[0,1]$ and some $\eps>0.$ This means that $m_{\vec{\beta}}$
is not concentrated near any one point $x\in[0,1]$
and, therefore,
$$
\Bigl(\int q^{p_1} dm_{\vec{\beta}}(q)\Bigr)^{1/p_1} \leq 
\Bigl(\int q^{p_2} dm_{\vec{\beta}}(q)\Bigr)^{1/p_2} - \delta
$$
for some $\delta>0.$
By (\ref{derapprox}), for large enough $N,$
$$
\bigl(\e\la R_{1,2}^{p_1}\ra\bigr)^{1/p_1}
\leq
\bigl(\e\la R_{1,2}^{p_2}\ra\bigr)^{1/p_2}
-\frac{\delta}{2}
$$
which means that the Gibbs measure can not concentrate near one point and,
therefore,
\begin{equation}
\e\bigl\la (R_{1,2} - \e \la R_{1,2}\ra)^2\bigr\ra\geq \delta'>0.
\label{RSBii}
\end{equation}
\qed

Even though these examples strengthen and generalize the result of Pastur and Shcherbina
in \cite{PS}, unfortunately, the argument used above does not apply to the most interesting case of
the classical $2$-spin model with external field, $\beta_2\not=0, h\not =0,$ 
and it is not clear how to prove (\ref{RSBii}) in that case.

\textbf{Acknowledgments.} The author would like to thank the referees for
many helpful comments and suggestions that lead to the improvement of the paper.

\end{document}

%% file: invlat.tex
\font\tencmmib=cmmib10 \skewchar\tencmmib '60
\newfam\cmmibfam
\textfont\cmmibfam=\tencmmib

\def\bbox{\quad\hbox{\vrule \vbox{\hrule \vskip2pt \hbox{\hskip2pt
\vbox{\hsize=1pt}\hskip2pt} \vskip2pt\hrule}\vrule}}
\def\lessim{\ \lower4pt\hbox{$
\buildrel{\displaystyle <}\over\sim$}\ }
\def\gessim{\ \lower4pt\hbox{$\buildrel{\displaystyle >}
\over\sim$}\ }

\def\eps{\varepsilon}

\def\go0{\to 0}

\def\la{\langle}

\def\leftitem#1{\item{\hbox to\parindent{\enspace#1\hfill}}}

\def\qed{\hfill\break\rightline{$\bbox$}}

\def\ra{\rangle}

\def\sg{\sigma}

\def\sg2{\sigma^2}

\def\__{_{\infty}}